\documentclass[12pt]{amsart}
\usepackage{amsmath, amsthm, amscd, amsfonts}
\usepackage{mathrsfs}
\usepackage{textcomp}
\usepackage{tikz-cd}
\usepackage[matrix,arrow]{xy}

\newtheorem{theorem}{Theorem}[section]
\newtheorem{lemma}[theorem]{Lemma}

\newtheorem{cor}[theorem]{Corollary}
\newtheorem{prop}[theorem]{Proposition}
\theoremstyle{definition}
\newtheorem{definition}[theorem]{Definition}
\theoremstyle{definition}
\newtheorem{remark}[theorem]{Remark}
\theoremstyle{definition}

\theoremstyle{definition}

\newcommand{\id}{\mbox{$\mathfrak{I}$}}
\newcommand{\co}{\mbox{$\mathscr{C}$}}
\newcommand{\rnd}{\mbox{$\mathfrak{R}$}}
\newcommand{\kel}{\mbox{$ \prec\ $}}
\newcommand{\sq}{\mbox{$ \sqsubseteq\ $}}
\newcommand\set[1]{\mbox{$\{#1\}$}}
\newcommand\cat[1]{\mbox{$\mathbf{#1}$}}
\newcommand\ds[1]{\mbox{$\displaystyle{#1}$}}
\newcommand\map[3]{\mbox{$#1:#2\rightarrow #3$}}
\newcommand\inv[1]{\mbox{$#1^{-1}$}}
\newcommand\drct[1]{\mbox{$\mathscr{#1}$}}
\newcommand\clos[1]{\mbox{$\overline{#1}$}}
\newcommand\comp[2]{\mbox{$#1\cdot #2$}}

\newcommand\astcomp[2]{\mbox{$#1\ast #2$}}
\newcommand\join[2]{\mbox{$#1\vee #2$}}
\newcommand\meet[2]{\mbox{$#1\wedge #2$}}
\newcommand\powcat[2]{\mbox{$#1^{\mathbb{#2}}$}}
\newcommand\kleicat[2]{\mbox{$#1_{\mathbb{#2}}$}}

\begin{document}

\title{A correspondence between proximity homomorphisms and certain frame maps via a comonad}

\author{Ando Razafindrakoto $^1$}

\address{$^{1}$ Department of Mathematics and Applied Mathematics, Uiversity of the Western Cape, Bellville 7535, Cape Town, South Africa.}

\email{arazafindrakoto@uwc.ac.za}

\subjclass[2010]{06D22, 18A40, 18C15, 18C20, 18B35, 54B30, 54E05, 54D35}

\keywords{Proximity, proximity homomorphism, round ideal, stably compact frame, compactification, comonad, coalgebras, Kleisli composition.}

\maketitle
{\centering\footnotesize {\em To Themba Dube on his $65^{\text{th}}$ Birthday}.\par}

\begin{abstract}
We exhibit the proximity frames and proximity homomorphisms as a Kleisli category of a comonad whose underlying functor takes a proximity frame  to its frame of round ideals. This construction is known in the literature as {\em stable compactification} (\cite{BezHar2}). We show that the frame of round ideals naturally carries with it two proximities of interest from which two comonads are induced. 
\end{abstract}

\section{Introduction}

The interest in proximities or strong inclusions for topological spaces lies in the fact that they allow a description of compactifications through binary relations on the powerset lattice. To our knowledge, this particular approach was pioneered by Smirnov in \cite{Smi}. Early studies on proximities include the work of V. A. Efremovi\v{c} (\cite{Efr}) and C. H. Dowker (\cite{Dow}), and their use have become important in understanding the property of compactness and total boundedness in the theory of (quasi-) uniform spaces. The transfer of this concept to the pointfree setting - hence replacing the powerset lattice with the lattice of open sets - has been investigated by A. Pultr and J. Picado in \cite{PicPul2}, J. Frith (\cite{Fri}) - who cites the work of S. A. Naimpally and B. D. Warrack (\cite{NaiWar}) among others, B. Banaschewski (\cite{Ban90}) and de Vries (\cite{deV}). The study of compactifications of frames or locales has led to two noticeable constructions: the use of {\em strong inclusions} on frames by B. Banaschewski (\cite{Ban90}) and the introduction of {\em proximities} on Boolean Algebras by de Vries (\cite{deV}). There appear two possible candidates for the morphisms to be considered between frames with such structures: one is that of a frame homomorphism that preserves the proximity relation and the relation as in \cite{Fri}, and another one is that of a subadditive meet-semilattice homomorphism that preserves the bottom element, as used by G. Bezhanishvili and J. Harding in \cite{BezHar,BezHar2}. The salient difference is that the latter may fail to be a frame map.

The present paper provides a convenient bridge between the two types of morphisms. Thus each {\em proximity homomorphism} \map{f}{L}{M} in the sense of \cite{BezHar,BezHar2} is shown to be equivalent to a frame homomorphism \map{\psi}{\rnd L}{M} that preserves proximities. Here $\rnd L$ is the frame of round ideals on $L$ endowed with the way below relation. This equivalence is presented on a categorical level and represents the frames with proximities and proximity homomorphisms as the Kleisli category of an idempotent comonad whose underlying functor is given by \rnd. The description of \rnd\ as a comonad also provides a functorial basis for the {\em stable compactification} of {\em proximity frames} introduced by G. Bezhanishvili and J. Harding (\cite{BezHar2,BezHar3}). 

The main results of the paper are in Section 3 and Section 4. While the bijection between the proximity morphisms \map{f}{L}{M} and special frame homomorphisms \map{\psi}{\rnd L}{M}, as mentioned in the previous paragraph, is shown in Section 3, it is in Section 4 that this correspondence is proved to be natural in the category theoretic sense. In particular, the functoriality of \rnd, as well as the study of the natural transformations that make it a comonad, are discussed in Section 4. We further show in this section that each $\rnd L$ carries with it a proximity that is maximally compatible with the proximity on $L$, and that this yields a comonad (generally non-idempotent) admitting \rnd\ as a sub-comonad. Generalities concerning proximity frames and monads are introduced in Section 2 in order to articulate those results in a convenient manner. The last section discusses classical and familiar compactifications in connection with the stable compactification. We wish to express the fact that the consideration of the relationship between $\rnd L$ and $\rnd\rnd L$, which is central to the paper, and which may be confusing on a superficial level, has been inspired by Lawvere's article \cite{Law} which uses a ``cylinder'' $\xymatrix{B \ar@<-.5ex>[r] \ar@<.5ex>[r] & I\times B\ar[r] & B}$ as a model to ``capture'' the dialectical principle of {\em unity and identity of opposites}. Where this may become relevant in Topology is when a pointed endofunctor $R$ with a universal property fails to be idempotent (\cite{Raz}). In such a case, the action of $R$ on the unit $1\to R$ still provides some useful information on $R$ and this is the case for a monad where we have such a cylindric model. 

\section{Preliminaries}

\paragraph*{\textbf{Proximities}.} A frame\footnote{We refer to \cite{Joh} and \cite{PicPul} for a general background.} $L$ is a complete lattice with the equational identity
\begin{center}
$\meet{a}{\ds{\bigvee_{i\in I}} b_i}=\ds{\bigvee_{i\in I}}\left(\meet{a}{b_i}\right)$
\end{center}
for all $a\in L$ and for any index set $I$. Frame homomorphisms or frame maps are functions \map{f}{L}{M} that preserve finite meets (including the top element $1$) and arbitrary joins (including the bottom element $0$). Frames and their homomorphisms form a category denoted by \cat{Frm}. An example of a frame is the lattice of open sets $\mathcal{O}X$ of a topological space $X$. In this case, any continuous function \map{f}{X}{Y} gives rise to a frame homomorphism \map{\inv{f}}{\mathcal{O}Y}{\mathcal{O}X} defined as the restriction of the inverse image map on the open sets. The abstract spaces for which the frames constitute formal open sets are called {\em locales}. They are the objects of the category $\cat{Loc}=\cat{Frm^{op}}$. Thus we have a functor \map{\mathcal{O}}{\cat{Top}}{\cat{Loc}}. For each frame $L$, one can assign the space $\Sigma L=\set{\map{f}{L}{2}\ |\ f\text{ is a frame map}}$ with the topology \set{\Sigma_a\ |\ a\in L}, where $\Sigma_a=\set{\map{f}{L}{2}\ |\ f(a)=1}$ consists of all points contained in $a\in L$. This induces a functor \map{\Sigma}{\cat{Loc}}{\cat{Top}} by assigning $\Sigma(f)(p)=pf$ for any frame map \map{f}{L}{M} and any point \map{p}{M}{2}. This gives an adjunction $\mathcal{O}\dashv \Sigma$.

\begin{definition}
\label{defn: definition of a proximity}
\cite{BezHar} A proximity \kel\ on a frame $L$ is a binary relation satisfying the following:
\begin{enumerate}
\item \kel\ is finer than $\leq$ and is a sublattice of $L\times L$;
\item If $a\leq b\kel c\leq d$, then $a\kel d$;
\item \kel\ is interpolative;
\item $a=\bigvee\set{b\ |\ b\kel a}$ for any $a\in L$.
\end{enumerate}
\end{definition}
The pair $(L,\prec)$ is called a {\em proximity frame}. Note that the difference with this definition and that of a strong inclusion (\cite{Ban90}) is in the following two axioms:
\begin{itemize}
\item If $a\kel b$, then $a$ is ``well inside'' $b$, that is, $\join{a^*}{b}=1$;
\item  If $a\kel b$, then $b^*\kel a^*$.
\end{itemize}
Thus a proximity frame is not necessarily regular. One such example is given by stably continuous frames (\cite{BanBru}). These are the frames endowed with the relation $\ll$ defined by
\begin{center}
$a\ll b$ if for any arbitrary $S\subseteq L$, when $b\leq \bigvee S$ then $a\leq\bigvee F$ where $F$ is finite and $F\subseteq S$.
\end{center}
The pair $(L,\ll)$ is a proximity frame that is not regular. The frame maps that preserve $\ll$ are called {\em proper} or {\em perfect}. This gives the category $\cat{StKFrm}$. There is a full embedding from \cat{StKFrm} to \cat{PrFrm} (\cite[Proposition 4.2]{BezHar}) and we directly use the order $\ll$ when we refer to those proximity frames that are in the range of this embedding.

\begin{definition}
\cite{BezHar} A {\em proximity homomorphism} from $(L,\prec)$ to $(M,\prec)$ is a meet-semilattice homomorphism \map{f}{L}{M} satisfying:
\begin{enumerate}
\item $f(0)=0$;
\item If $a_1\kel b_1$ and $a_2\kel b_2$, then $f(\join{a_1}{a_2})\kel\join{f(b_1)}{f(b_2)}$;
\item $f(a)=\bigvee\set{f(b)\ |\ b\kel a}$.
\end{enumerate}
\end{definition}
Proximity frames and proximity homomorphisms form a category \cat{PrFrm} with the composition given by
\begin{center}
$(\astcomp{g}{f})(a)=\bigvee\set{g(f(b))\ |\ b\kel a}$.
\end{center}
As shown in \cite[Example 3.5]{BezHar}, proximity homomorphisms do not need to be frame homomorphisms\footnote{In particular, they do not necessarily preserve finite joins.}. In general, we have $(\astcomp{g}{f})(a)\leq (gf)(a)$. We are naturally interested in the category \cat{C} of proximity frames and frame maps that preserve proximities.

\begin{lemma}
\label{lem: When prox morphisms composition coincide with usual composition}
(\cite[Lemma 3.7]{BezHar}) If $g$ preserves joins, then $\astcomp{g}{f}=gf$.
\end{lemma}

\paragraph*{Monads.} A {\em monad} $\mathbb{T}$ on a category \cat{C} is a triple $(T,\mu,\eta)$, where $T$ is an endofunctor on \cat{C}, \map{\mu}{TT}{T} and \map{\eta}{1}{T} are natural transformations satisfying the identities 
\begin{center}
$\comp{\mu}{T\mu}=\comp{\mu}{\mu T}$ and $\comp{\mu}{\eta T}=\comp{\mu}{T\eta}=1_T$.
\end{center}
A $\mathbb{T}$-algebra (or an {\em Eilenberg-Moore algebra}) is a pair $(X,a)$, where $X\in\cat{C}$ and \map{a}{TX}{X} a morphism such that 
\begin{center}
$\comp{a}{Ta}=\comp{a}{\mu_X}$ and $\comp{a}{\eta_X}=1_X$.
\end{center}
Note that $(TX,\mu_X)$ is the free $\mathbb{T}$-algebra over $X$. If $(X,a)$ and $(Y,b)$ are $\mathbb{T}$-algebras, then a $\mathbb{T}$-algebra homomorphism \map{f}{(X,a)}{(Y,b)} is a morphism \map{f}{X}{Y} in \cat{C} such that $\comp{f}{a}=b\cdot Tf$. The category of $\mathbb{T}$-algebras and $\mathbb{T}$-algebra homomorphisms are denoted by \powcat{\cat{C}}{T}. The forgetful functor $\map{\powcat{G}{T}}{\powcat{\cat{C}}{T}}{\cat{C}}:(X,a)\mapsto X$ admits a left adjoint $\map{\powcat{F}{T}}{\cat{C}}{\powcat{\cat{C}}{T}}:X\mapsto (TX,\mu_X), f\mapsto Tf$. The unit of this adjunction is given by $\eta_X:X\to \powcat{G}{T}\powcat{F}{T}X$ and the co-unit is provided by structure maps $\varepsilon_{(X,a)}:\powcat{F}{T}\powcat{G}{T}(X,a)\to(X,a)$. 

The {\em Kleisli category} \kleicat{\cat{C}}{T} associated to $\mathbb{T}$ consists of the same objects as those of \cat{C}, and arrows $X\rightharpoonup Y$ are morphisms $X\to TY$ in \cat{C}. The composition $g\bullet f$ of \map{f}{X}{TY} and \map{g}{Y}{TZ} is given by the morphism $\mu_Z\cdot Tg\cdot f$. The functor \map{\kleicat{F}{T}}{\cat{C}}{\kleicat{\cat{C}}{T}} defined by $\kleicat{F}{T}(X)=X$ and $\kleicat{F}{T}(\map{f}{X}{Y})=Tf\cdot \eta_X$ admits a right adjoint \map{\kleicat{G}{T}}{\kleicat{\cat{C}}{T}}{\cat{C}} which takes \map{f}{X}{TY} to \map{\mu_Y\cdot Tf}{TX}{TY}. The dual notion of a monad is called {\em comonad} and that of an algebra is {\em coalgebra}; we will however use the same term {\em Kleisli} to refer to the ``co-Kleisli category'' of a comonad. For further generalities on monads, we refer the reader to the texts \cite{Mac}, \cite{Bor2} and \cite{Low}.

\section{The frame of round ideals}
The feature of a proximity morphism \map{f}{L}{M}, namely the identity $f(a)=\bigvee\set{f(b)\ |\ b\kel a}$, allows us to decompose it as follows
\begin{center}
$a\mapsto \set{b\ |\ b\kel a}\mapsto\set{f(b)\ |\ b\kel a}\mapsto \bigvee\set{f(b)\ |\ b\kel a}$.
\end{center}
This necessitates that we concentrate on round ideals.
\begin{definition}
A subset $J\subseteq L$ is called an {\em ideal} if $J$ is a downset and closed under the formation of finite joins. The set of all ideals on a frame $L$ is denoted by $\id L$. $\id L$ is a frame with the following operations: meets are given by set-intersections and joins given by
\begin{center}
$\bigvee\mathscr{J}=\bigcup\set{I_1\vee I_2\vee\dots\vee I_n\ |\ I_1, I_2, \dots, I_n\in\mathscr{J}\text{ and }n\in \mathbb{N}}$,
\end{center}
where $I_1\vee I_2\vee\dots\vee I_n=\set{i_1\vee i_2\vee\dots\vee i_n\ |\ i_k\in I_k\text{ for }1\leq k\leq n}$. Thus, if $\mathscr{J}$ is a directed set, then $\bigvee\mathscr{J}=\bigcup\mathscr{J}$.
\end{definition}

\begin{definition}
On a proximity frame $L$, an ideal $I$ is said to be {\em round} if for any $a\in I$, there is $b\in I$ such that $a\kel b$. The collection of round ideals in $L$ is denoted by $\rnd L$. 
\end{definition}
\begin{lemma}
\label{lem: round ideals form a frame}
$\rnd L$ is a subframe of the frame of ideals $\id L$. The composition \map{\varsigma_L}{\rnd L}{L} of the frame inclusion $m:\rnd L\to \id L$ and the join map \map{\bigvee}{\id L}{L} is then a frame homomorphism.
\end{lemma}
Lemma \ref{lem: round ideals form a frame} can be directly seen in \cite[Proposition 4.6 (1)]{BezHar}. In particular we mention that $\rnd L$ is stably compact (\cite[Proposition 4.8]{BezHar}) with $I\ll J$ if and only if $m(I)\ll m(J)$.
\begin{prop}
\label{prop: semi-lattice hom to semi-lattice hom}
If \map{f}{L}{M} is a meet-semilattice homomorphism between two proximity frames, then the function \map{\rnd f}{\rnd L}{\rnd M} given by
\begin{center}
$\rnd f(I)=\set{a\ |\ a\kel f(b)\text{ for some }b\in I}$,
\end{center}
is a meet-semilattice homomorphism.
\end{prop}
\begin{proof}
Let us first show that $\rnd f$ is well-defined. Notice that $\rnd f(I)$ is a downset. Now, let us take $a,b\in \rnd f(I)$. For some $c,d\in I$,\ $a\kel f(c)$ and $b\kel f(d)$. Since $f$ is monotone, we have
\begin{center}
$\join{a}{b}\kel \join{f(c)}{f(d)}\leq f(\join{c}{d})$,
\end{center}
showing that $\join{a}{b}\in\rnd f(I)$. Now let $a\in \rnd f(I)$ and $b\in I$ with $a\kel f(b)$. Since \kel\ is interpolative, there is $c\in M$ such that $a\kel c\kel f(b)$, showing that $\rnd f(I)$ is round. Thus $\rnd f$ is well-defined.

Now let us show that $\rnd f$ preserves finite meets. Let $a\in\rnd f(I)\cap \rnd f(J)$. We have $a\kel\meet{f(c)}{f(d)}=f(\meet{c}{d})$ for some $c\in I$ and $d\in J$. Hence $a\in\rnd f(I\cap J)$. It is straightforward to check that $\rnd f(L)=M$.
\end{proof}

\begin{remark}
\label{rem: remarks about joins}
As can be noticed above, if $f$ preserves \kel, then $f(I)\subseteq\rnd f(I)\subseteq\id f(I)$. However, even if this is not true, the difference cannot be detected by joins:
\begin{center}
$\ds{\bigvee} f(I)=\ds{\bigvee}\set{a\ |\ a\leq f(b)\text{ for some }b\in I}=\ds{\bigvee}\rnd f(I)$.
\end{center}
\end{remark}

\begin{prop}
\label{prop: semi-lattice to frame hom}
If $L$ and $M$ are two frames and \map{f}{L}{M} is a meet-semilattice homomorphism that preserves \kel, then the function \map{\rnd f}{\rnd L}{\rnd M} is a proper frame homomorphism.
\end{prop}

\begin{proof}
Let us show that $\rnd f$ preserves finite and directed joins. Let $a\in\rnd f(\join{I}{J})$ and let $c\in I$ and $d\in J$ such that $a\kel f(\join{c}{d})$. There are $c'\in I$ and $d'\in J$ such that $c\kel c'$ and $d\kel d'$. We then have $a\kel f(\join{c}{d})\kel\join{f(c')}{f(d')}$ and so $a\in\join{\rnd f(I)}{\rnd f(J)}$. On the other hand, $\rnd(f)(\set{0})=\set{0}$.

Let \drct{J} be a directed subset and let $a\in \rnd f(\bigcup\drct{J})$. There are $K\in\drct{J}$ and $b\in K$ such that $a\kel f(b)$. Thus 
\begin{center}
$a\in\rnd f(K)\subseteq\bigcup\set{\rnd f(K)\ |\ K\in\drct{J}}=\bigvee\set{\rnd f(K)\ |\ K\in\drct{J}}$.
\end{center}

Now suppose that $I\ll J$ and let \drct{J} be a directed subset such that $\rnd f(J)\subseteq\bigcup\drct{J}$. There is $a\in J$ such that $I\subseteq \downarrow a$. Since $f$ preserves \kel, we have $f(a)\in \rnd f(J)\subseteq\bigcup\drct{J}$. There is $K\in\drct{J}$ such that for all $i\in I$,\ $f(i)\leq f(a)\in K$. Therefore $\rnd f(I)\subseteq K$. This shows that $\rnd f(I)\ll\rnd f(J)$.
\end{proof}
\begin{cor}
\label{cor: proxi map to frame map}
For any given proximity homomorphism \map{f}{L}{M}, the function \map{\rnd f}{\rnd L}{\rnd M} is a proper frame homomorphism.
\end{cor}

Now, consider the map \map{\kappa_L}{L}{\rnd L} where $\kappa_L(a)=\set{b\ |\ b\kel a}$ for each frame $L$. $\kappa_L(a)$ is the largest round ideal such that $\varsigma_L(\kappa_L(a))=a$ and is in fact the right adjoint of $\varsigma_L$. (See \cite{BezHar} and \cite[Theorem 4.20]{BezHar2}). From \cite[Proposition 4.6.]{BezHar}, we know that $\kappa$ is a proximity homomorphism that takes \kel\ to $\ll$.  On the other hand, the join \map{\varsigma_L}{\rnd L}{L} is a frame homomorphism that takes $\ll$ to the proximity $\kel$. We shall only use the notations $\varsigma$ and $\kappa$ when there is no risk of confusion.

\begin{lemma}
\label{lem: kappa natural up to joins}
For any meet-semilattice homomorphism \map{f}{L}{M} between proximity frames satisfying the property $f(a)=\bigvee\set{f(b)\ |\ b\kel a}$, we have $f=\varsigma_M\rnd f\kappa_L$.
\end{lemma}
\begin{proof}
Let $a\in L$. We have
\begin{align*}
\varsigma_M(\rnd f(\kappa_L(a))) &=\bigvee\set{c\ |\ c\kel f(b)\text{ for some } b\kel a}\\
&=\bigvee\set{\ds{\bigvee}\set{c\ |\ c\kel f(b)}\ |\ \text{ for some } b\kel a}\\
&= \bigvee\set{f(b)\ |\ \text{ for some } b\kel a}\\
&= f(a).
\end{align*}
As $a$ is arbitrary, the result follows.
\end{proof}

\begin{remark}
\label{rem: expression of composition of proximity homomorphism}
Note that we can write $\astcomp{g}{f}=\comp{\varsigma}{\rnd(gf)}\cdot\kappa$. This is only equal to $gf$ when the latter is a proximity homomorphism.
\end{remark}

\begin{prop}
\label{prop: stably compact proximity frame}
If $L$ is a stably compact proximity frame, then $\varsigma$ admits a left adjoint \map{\alpha}{L}{\rnd L} such that $\alpha(a)\subseteq\kappa(a)$ for all $a\in L$. 
\end{prop}
\begin{proof}
Define $\alpha(a)=\set{b\ |\ b\ll a}$. We have $\alpha(a)\subseteq\kappa(a)$ because $a=\varsigma(\kappa(a))$. This shows that $\alpha(a)$ is round. Since $L$ is stably compact, we have $a=\varsigma(\alpha(a))$. On the other hand, if $a\leq\varsigma(J)$, then $\alpha(a)\subseteq J$ and so $\alpha(\varsigma(J))\subseteq J$. 
\end{proof}
The above result shows that if we consider $\rnd L$ as an object of \cat{PrFrm} or \cat{C} by taking $\ll$ as a proximity, then $\rnd\rnd L\cong\rnd L$ both in \cat{PrFrm} (\cite[Proposition 4.9]{BezHar}) and \cat{C}. This is due to the fact that the left and right adjoint of $\varsigma_{\mathfrak{R} L}$ coincide. On the other hand if we take $\id L$, the frame of $\leq$-round ideals (\cite{BezHar2}), and if we consider the natural proximity $\subseteq$ on $\id L$,  $(\id(\id L,\subseteq),\subseteq)$ is not isomorphic to $(\id L,\subseteq)$.

\begin{lemma}
\label{lem: theta well-defined}
For each proximity homomorphism \map{f}{L}{M}, the frame homomorphism \map{\varsigma_M\rnd f}{\rnd L}{M} takes $\ll$ to \kel.
\end{lemma}
\begin{lemma}
\label{lem: rho well-defined}
For each frame homomorphism $\map{\psi}{\rnd L}{M}$ that takes $\ll$ to \kel, the composition $f=\psi\kappa_L$ is a proximity homomorphism.
\end{lemma}
\begin{proof}
Straightforward verification (\cite{BezHar}).
\end{proof}
We have shown that the assignments 
\begin{align*}
  \theta_{L,M}:\ \cat{PrFrm}(L;M)& \longrightarrow \cat{C}(\rnd L;M) \\
  f & \longmapsto \varsigma_M\rnd f
\end{align*}
and
\begin{align*}
  \rho_{L,M}:\ \cat{C}(\rnd L;M) & \longrightarrow \cat{PrFrm}(L;M) \\
  \varphi & \longmapsto \varphi\kappa_L
\end{align*}
are well-defined. They are in fact inverse to each other.

\begin{lemma}
\label{lem: from frame map back to frame map}
If \map{\psi}{\rnd L}{M} is a frame homomorphism, then $\psi=\varsigma_M\rnd (\psi\kappa_L)$.
\end{lemma}
\begin{proof}
Note that for any round ideal $I$,\ $\psi(I)=\bigvee\set{\psi\kappa_L(b)\ |\ b\in I}$. We have that
\begin{align*}
\varsigma_M\rnd (\psi\kappa_L)(I) &= \bigvee\set{b\ |\ b\kel \psi\kappa_L(c)\text{ for some }c\in I}\\
&= \bigvee\set{\ds{\bigvee}\set{b\ |\ b\kel \psi\kappa_L(c)}\ |\ c\in I}\\
&= \bigvee\set{\psi\kappa_L(c)\ |\ c\in I}\\
&= \psi(I),
\end{align*}
showing the equality.
\end{proof}
\begin{theorem}
\label{thm: one-one correspondence between morphisms}
$\comp{\theta_{L,M}}{\rho_{L,M}}=1$ and $\comp{\rho_{L,M}}{\theta_{L,M}}=1$.
\end{theorem}
\begin{proof}
Lemma \ref{lem: kappa natural up to joins} and Lemma \ref{lem: from frame map back to frame map}.
\end{proof}

The problem of determining whether the bijection in Theorem \ref{thm: one-one correspondence between morphisms} is natural in the category theoretic sense requires an understanding of either $\kappa_L$ or $\varsigma_M$. G. Bezhanishvili and J. Harding have shown in \cite{BezHar} that $\astcomp{\kappa_L}{\varsigma_L}=1_{\mathfrak{R} L}$ and $\astcomp{\varsigma_L}{\kappa_L}=1_L$, thereby showing that $\rnd L\cong L$ in \cat{PrFrm}. In Lemma \ref{lem: from frame map back to frame map}, we actually have $\psi\kappa=\varsigma\rnd(\psi\kappa)\kappa$ as in Lemma \ref{lem: kappa natural up to joins}. In light of the previous identities that imply that $\kappa$ is an epimorphism in \cat{PrFrm}, this reduces to the result $\psi=\varsigma\rnd(\psi\kappa)$.

We shall show in the next section that this isomorphism is the concrete realisation of the representation of proximity frames as coalgebras.

\section{Functorial constructions}
We first consider the structure $(\rnd L,\ll)$ - that we shall simply denote by $\rnd L$. The structure $\rnd L$ is known as {\em the stable compactification} of $(L,\prec)$ and it is shown in \cite{BezHar2} that this generalises the results of Banaschewski (\cite{Ban90}) on the equivalence between frame compactifications and strong inclusions to the context of proximity frames. 

\subsection{Stable compactification}

\begin{lemma}
\label{lem: m is a natural transformation}
For a proximity morphism \map{f}{L}{M}, we have $\comp{\id f}{m_L}=\comp{m_M}{\rnd f}$, where $m$ is the frame inclusion from Lemma \ref{lem: round ideals form a frame}.
\end{lemma}
\begin{proof}
Let $I$ be a round ideal. Since $f$ preserves $\prec$, $f(I)\subseteq \rnd f(I)\subseteq \id f(I)$, that is $(\comp{m_M}{\rnd f})(I)\subseteq(\comp{\id f}{m_L})(I)$. Now, if $a\leq f(c)$ for some $c\in I$, then $f(c)\kel f(d)$ for some $d\in I$. This shows that $(\comp{\id f}{m_L})(I)\subseteq(\comp{m_M}{\rnd f})(I)$.
\end{proof}
\begin{lemma}
\label{lem: rnd join is a nat transformation}
For a frame homomorphism \map{f}{L}{M} that preserves the relation $\prec$,\ $f\cdot\varsigma_L=\varsigma_M\cdot\rnd f$.
\end{lemma}
\begin{proof}
Consider the composition $\varsigma=\comp{\ds{\bigvee}}{m}$. Since $f$ is a frame homomorphism, we have $f\cdot\varsigma_L=\comp{f}{\comp{\ds{\bigvee}_L}{m_L}}=\comp{\ds{\bigvee}_M}{\comp{\id f}{m_L}}$. The result follows from Lemma \ref{lem: m is a natural transformation} above.
\end{proof}
\begin{cor}
\label{cor2a: proxi map to frame map}
\rnd\ is an endofunctor on \cat{C}
\end{cor}
\begin{proof}
For a composable pair $f,g$ in \cat{C}, it is clear that  $(\comp{\rnd g}{\rnd f})(I)\subseteq\rnd(gf)(I)$ for any round ideal $I$. Now, since $I$ is round, the reverse inclusion holds and for a proximity frame $L$, we have $\rnd(1_L)(I)=1_{\mathfrak{R} L}(I)$.
\end{proof}

\begin{remark}
\label{rem: on proxi map to frame map}
If $g\ast f=gf$, then $\rnd(g\ast f)=\rnd(gf)=\rnd g\rnd f$. The last equality is in fact always true for meet-semilattice homomorphisms that preserve $\prec$.
\end{remark}
\begin{cor}
\label{cor2: s and c join is a nat transformation}
$\varsigma:\rnd\to 1$ is a natural transformation.
\end{cor}

We now consider the assignments $r_L=\rnd(\kappa_L)$.
We have 
\begin{align*}
\rnd(\kappa_L)(I) &= \set{K\in \rnd L\ |\ K\ll\kappa_L(a)\text{ for some }a\in I}\\
 &= \set{K\in \rnd L\ |\ \varsigma_L(K)\kel a\text{ for some }a\in I}\\
 &= \set{K\in \rnd L\ |\ \varsigma_L(K)\in I} \text{ (Since }I\text{ is round)}\\
 &= \set{K\in \rnd L\ |\ K\ll I} \text{ (I}\text{ is round)}
\end{align*}
By Proposition \ref{prop: stably compact proximity frame}, \map{r_L}{\rnd L}{\rnd\rnd L} is an isomorphism. $r$ is in fact a natural isomorphism:
\begin{prop}
\label{prop: comultiplication is a nat transformation}
For each proximity morphism \map{f}{L}{M}, we have\ $\rnd\rnd f\cdot r_L=r_M\cdot \rnd f$.
\end{prop}
\begin{proof}
Note that for each $I\in\rnd L$,
\begin{center}
$r_M(\rnd f(I)) = \set{ J\ |\ \varsigma_M(J)\kel f(a)\text{ for some }a\in I}$. 
\end{center}
Thus we have
\begin{align*}
\rnd\rnd f(r_L(I)) &= \set{K\ |\ K\ll \rnd f(J)\text{ for some }J\in r_L(I)}\\
&= \set{K\ |\ \varsigma_M(K)\in\rnd f(J)\text{ for some }J\in r_L(I)}\\
&= \set{K\ |\ \varsigma_M(K)\kel f(b)\text{ for some }b\in J \text{ and for some }J\ll I}\\
&\subseteq r_M(\rnd f(I)).
\end{align*}
The reverse inclusion in the last step holds since $I$ is a round ideal. 
\end{proof}

\begin{theorem}
\label{thm: R is a coreflector}
The triple $(\rnd, r,\varsigma)$ forms an idempotent comonad, i.e. $(\rnd,\varsigma)$ is essentially a coreflector from \cat{C} to \cat{StKFrm}.
\end{theorem}

\subsubsection{Kleisli category associated to the stable compactification}
The Kleisli category \kleicat{\cat{C}}{R} associated  to $(\rnd,r,\varsigma)$ is formed with the same objects as in \cat{C} but with morphisms $L\rightharpoonup M$ which are morphisms $\rnd L\to M$ in \cat{C}. The composition of \map{u}{\rnd L}{M} and \map{v}{\rnd M}{N} is then given by 
\begin{center}
${v}\bullet{u}=v\cdot \rnd u\cdot r_L$.
\end{center}
Let us define a functor \map{F}{\cat{PrFrm}}{\kleicat{\cat{C}}{R}}. Since each proximity homomorphism \map{f}{L}{M} can be expressed as $f=\varsigma_M\rnd (f)\kappa_L$, let $F(f)=\varsigma_M\rnd f=\theta_{L,M}(f)$. Beside the straightforward identity $F(1_L)=\varsigma_L=1_L:L\rightharpoonup L$, we need to show that
\begin{center}
$F(\astcomp{g}{f})= F(g)\bullet F(f)=\theta_{M,N}(g)\bullet\theta_{L,M}(f)$, 
\end{center}
for any $g,f$ in \cat{PrFrm}. This is given in the following proposition.
\begin{prop}
\label{prop: Kleisli identity}
Given $g,f$ in \cat{PrFrm}, we have 
\begin{center}
$\comp{\varsigma_N}{\rnd g}\cdot\rnd(\varsigma_M\cdot\rnd f)\cdot r_L=\varsigma_N\cdot\rnd(gf)$.
\end{center}
\end{prop}
\begin{proof}
Consider the diagram below where each rectangle commutes.
$$\xymatrix{
 & \rnd L\ar[rr]^{\mathfrak{R} f} && \rnd M\ar[rr]^{\mathfrak{R} g} && \rnd N\ar[r]^{\varsigma_N} & N\\ 
\rnd L\ar[r]_{r_L} & \rnd\rnd L\ar[u]^{\varsigma_{\mathfrak{R} L}} \ar[rr]_{\mathfrak{R}\mathfrak{R} f} &&  \rnd\rnd M\ar[u]^{\varsigma_{\mathfrak{R} M}} \ar[rr]_{\mathfrak{R}\mathfrak{R} g} \ar[d]_{\mathfrak{R}(\varsigma_M)} && \rnd\rnd N\ar[u]_{\varsigma_{\mathfrak{R} N}} \ar[d]^{\mathfrak{R}(\varsigma_N)}\\
 & && \rnd M\ar[rr]_{\mathfrak{R} g} && \rnd N\ar[r]_1 & \rnd N \ar[uu]_{\varsigma_N}
}$$
The rectangle at the bottom commutes from the facts that $\rnd(\varsigma)\cdot r=1$ and that $r_M$ is an isomorphism, hence an epimorphism. The outer diagram, with Remark \ref{rem: on proxi map to frame map}, gives the desired identity.
\end{proof}

\begin{theorem}
\label{thm: PrFrm as a Kleisli category}
$F$ is an isomorphism, i.e. \cat{PrFrm} appears as the Kleisli category associated to the comonad $(\rnd,r,\varsigma)$ on the category \cat{C}.
\end{theorem}
\begin{proof}
Theorem \ref{thm: one-one correspondence between morphisms}.
\end{proof}
We have the result (\cite{BezHar}) by G. Bezhanishvili and J. Harding  that \cat{Prfrm} is equivalent to \cat{StkFrm}. Indeed, since \rnd\ is idempotent, \kleicat{\cat{C}}{R} and \powcat{\cat{C}}{R} essentially coincide. Consequently, as the Kleisli category is equivalent to the subcategory of free coalgebras in general, \cat{PrFrm} is equivalent to the subcategory of coalgebras $\rnd L$.

\subsection{Maximal proximity associated to the stable compactification}
We consider the structure $\co L=(\rnd L,\sq)$ where \sq\ is defined below.
\begin{prop}
\label{prop: construction of a new proximity}
The relation \sq\ defined by $I\sq J$ if and only if $I\subseteq J$ and $I\ll \kappa(\varsigma(J))$ is a proximity relation on $\rnd L$. The relation \sq\ is maximal in a sense that $I\sq J$ if and only if $I\subseteq J$ and $\varsigma(I)\kel\varsigma(J)$.
\end{prop}
\begin{proof}
Clearly  $\set{0}\sq \set{0}$, $L\sq L$ and \sq\ is finer than $\subseteq$. If we have $K\subseteq I\sq J\subseteq L$, then $K\subseteq I$ and $I\ll \kappa(\varsigma(J))\ll \kappa(\varsigma(L))$ so that $K\sq L$. If $I,J\sq L$ and $I\sq M,N$, then $I\vee J\ll \kappa(\varsigma(L))$ and since $\kappa$ preserves finite meets $I\ll \kappa(\varsigma(M))\cap \kappa(\varsigma(N))=\kappa(\varsigma(M\cap N))$. Suppose now that $I\sq J$. Let $a$ such that $\varsigma(I)\kel a\kel\varsigma(J)$ and let $C=\kappa(a)\cap J$. We have $I\subseteq C\subseteq J$ and $I\ll\kappa(a)\ll\kappa(\varsigma(J))$. Therefore
\begin{center}
 $I\ll\kappa(\varsigma(\kappa(a)\cap J))=\kappa(a)$ and $C\subseteq\kappa(a)\ll\kappa(\varsigma(J))$.
\end{center}
This shows that $I\sq C\sq J$. Finally, since 
\begin{center}
$\bigvee\set{K\ |\ K\ll I}\subseteq\bigvee\set{K\ |\ K\ll \kappa(\varsigma(I))}$,
\end{center}
we have $I\subseteq\bigvee\set{K\ |\ K\sq I}\subseteq I$.
\end{proof}

\begin{prop}
\label{prop: Rf continuous with respect to sq}
For any frame homomorphism \map{f}{L}{M} that preserves $\prec$, $\rnd f$ preserves \sq.
\end{prop}
\begin{proof}
Suppose $I\ll \kappa_L(\varsigma_L(J))$. There is $a\in L$ such that $I\ll \kappa_L(a)\ll\kappa_L(\varsigma_L(J))$. We  then have $a=\varsigma_L(\kappa_L(a))\kel\varsigma_L(J)$. For any $j\in J$, $\meet{a}{j}\in J$ and $\meet{a}{j}=\bigvee\set{c\in J\ |\ c\kel\meet{a}{j}}$. On the other hand we have that $a=\bigvee\set{\meet{a}{j}\ |\ j\in J}$ and so $a=\bigvee\set{\bigvee\set{c\in J\ |\ c\kel\meet{a}{j}}\ |\ j\in J}=\bigvee\set{c\in J\ |\ c\kel a}$. Since $f$ preserves joins, 
\begin{center}
$f(a)=\bigvee\set{f(c)\ |\ c\in J\text{ and }c\kel a}=\varsigma_M\rnd f(\kappa_L(a)\cap J)\leq\varsigma_M\rnd f(J)$.
\end{center}
Since $f=\varsigma_M\rnd f\kappa_L$ we have $\varsigma_M\rnd f(I)\kel (\varsigma_M\rnd f\kappa_L)(a)\leq \varsigma_M\rnd f(J)$. This shows  $\rnd f(I)\ll\kappa_M\varsigma_M\rnd f(J)$.
\end{proof}
\begin{cor}
\label{cor2: proxi map to frame map}
\co\ is an endofunctor on \cat{C} and there is a natural bijection $\map{\beta}{\rnd}{\co}$.
\end{cor}
\begin{proof}
For any $\map{f}{L}{M}$ in \cat{C}, it has been shown that $\rnd f$ preserves both $\ll$ and \sq. For each proximity frame $L$, $\beta_L$ is defined as the inclusion $(\rnd L,\ll)\to(\rnd L, \sq)$. $\co f$ is then defined as the ``extension'' of $\rnd f$, i.e. $\co f\cdot \beta_L=\beta_M\cdot\rnd f$.
\end{proof}
\begin{cor}
\label{cor2: c join is a nat transformation}
There is a natural transformation $\varepsilon:\co\to 1$, where $\varepsilon\beta=\varsigma$.
\end{cor}
\begin{remark}
\label{rem: differentiating the elements of the two structures}
To differentiate the members of $\co L$ from those of $\rnd L$, we shall write $\clos{J}=\beta_L(J)$. The only exception where this will not apply is for $\drct{J}\in\co(\co L)$. $\co L$ is a stably compact frame and we have $\clos{I}\ll\clos{J}$ if and only if $I\ll J$. Also, we have $a\in I$ if and only if $a\in\clos{I}$. $\rnd L$ and $\co L$ are then isomorphic as frames, but not as proximity frames. 
\end{remark}
\begin{lemma}
\label{lem: sequence of adjunction from CL to L}
For a proximity frame $L$, $\varepsilon_L\dashv \beta_L\kappa_L$. If in addition $L$ is stably compact, then $\beta_L\alpha\dashv \varepsilon_L$, where \map{\alpha}{L}{\rnd L} is the monotone map from Proposition \ref{prop: stably compact proximity frame}.
\end{lemma}
\begin{proof}
By construction, $\comp{\varepsilon_L}{(\beta_L\kappa_L)}=1$ and $\comp{\varepsilon_L}{(\beta_L\alpha)}=1$. From the adjunction $\varsigma_L\dashv\kappa_L$, it follows $\beta_L\cdot(\kappa_L\varepsilon_L\beta_L)\geq\beta_L$. Since $\beta_L$ is surjective $\comp{(\beta_L\kappa_L)}{\varepsilon_L}\geq 1$. In the same way, $\comp{(\beta_L\alpha)}{\varepsilon_L}\leq 1$.
\end{proof}

Consider $c_L=\co(\beta_L\ast\kappa_L)=\co(\beta_L\kappa_L)$ (Lemma \ref{lem: When prox morphisms composition coincide with usual composition}):
\begin{align*}
\co(\beta_L\kappa_L)(\clos{I})  &= \co(\beta_L\kappa_L)\beta_L(I)\\
 &= \beta_{\mathscr{C}L}\rnd(\beta_L\kappa_L)(I)\\
 &= \beta_{\mathscr{C}L}\set{\clos{K}\in \co L\ |\ \clos{K}\ \sq\beta_L\kappa_L(a)\text{ for some }a\in I}\\
 &= \beta_{\mathscr{C}L}\set{\clos{K}\in \co L\ |\ \varepsilon_L(\clos{K})\kel a\text{ for some }a\in I}\\
 &= \beta_{\mathscr{C}L}\set{\clos{K}\in \co L\ |\ \varepsilon_L(\clos{K})\in I} \text{ (Since }I\text{ is round)}
 \end{align*}
\begin{prop}
\label{prop: c is a morphism in C}
Each $c_L=\beta_{\mathscr{C}L}\rnd(\beta_L\kappa_L)$ is a morphism in \cat{C} and we have $c_L\dashv\varepsilon_{\mathscr{C}L}\dashv\beta_{\mathscr{C}L}\kappa_{\mathscr{C}L}$.
\end{prop}
\begin{proof}
Note that 
\begin{center}
$\set{\clos{K}\in \co L\ |\ \varepsilon_L(\clos{K})\in I}=\set{\clos{K}\in \co L\ |\ K\ll I}=\set{\clos{K}\in \co L\ |\ \clos{K}\ll \clos{I}}$.
\end{center}
Therefore $c_L=\beta_{\mathscr{C}L}\alpha$ where \map{\alpha}{\co L}{\rnd\co L}. The result follows from Lemma \ref{lem: sequence of adjunction from CL to L}.
\end{proof}
\begin{prop}
\label{prop: c is a natural transformation}
The morphisms $c_L$ define a natural transformation \map{c}{\co}{\co\co}.
\end{prop}
\begin{proof}
Let \map{f}{L}{M} be a morphism in \cat{C} and let $\clos{I}\in\co L$. We have
\begin{align*}
\co\co f(c_L(\clos{I})) &=\beta_{\mathscr{C}M}\set{\clos{K}\in\co L\ |\ \clos{K}\ \sq \co f(\clos{J})\text{ with }\clos{J}\in c_L(\clos{I}) }\\
&=\beta_{\mathscr{C}M}\set{\clos{K}\in\co L\ |\ \clos{K}\ \sq \beta_M(\rnd f(J))\text{ with }\varepsilon_L(\clos{J})\in \clos{I}}\\
&=\beta_{\mathscr{C}M}\set{\clos{K}\in\co L\ |\ \varepsilon_M(\clos{K})\kel\varsigma_L(\rnd f(J))\text{ with }\varepsilon_L(\clos{J})\in \clos{I}}\\
&=\beta_{\mathscr{C}M}\set{\clos{K}\in\co L\ |\ \varepsilon_M(\clos{K})\kel f(\varsigma_L(J))\text{ with }\varepsilon_L(\clos{J})\in \clos{I}}.\\
&\text{ (Lemma \ref{lem: rnd join is a nat transformation})}
\end{align*} 
On the other hand
\begin{align*}
c_M(\co f(\clos{I})) &=\beta_{\mathscr{C}M}\set{\clos{K}\in\co L\ |\ \varepsilon_M(\clos{K})\in \co f(\clos{I})}\\
&=\beta_{\mathscr{C}M}\set{\clos{K}\in\co L\ |\ \varepsilon_M(\clos{K})\in \co f(\beta_L(I))}\\
&=\beta_{\mathscr{C}M}\set{\clos{K}\in\co L\ |\ \varepsilon_M(\clos{K})\in \beta_M(\rnd f(I))}\\
&=\beta_{\mathscr{C}M}\set{\clos{K}\in\co L\ |\ \varepsilon_M(\clos{K})\kel f(i),\ i\in I}.
\end{align*}
The two are equal by considering $J=\kappa_L(i)$.
\end{proof}
\begin{lemma}
\label{lem: convexity and comultiplication}
For each $\drct{J}\in\co\co L$ and $\clos{I}\in\co L$,\ $\varepsilon_L(\varepsilon_{\mathscr{C}L}(\drct{J}))\in \clos{I}$ if and only if for some $\clos{J}\in \co L$, $\varepsilon_{\mathscr{C}L}(\drct{J})\ \sq \clos{J}$ and $\varepsilon_L(\clos{J})\in \clos{I}$.
\end{lemma}
\begin{proof}
If $\varepsilon_L(\varepsilon_{\mathscr{C}L}(\drct{J}))\kel b$ for some $b\in \clos{I}$, then we take $\clos{J}=\beta_L\kappa_L(b)$. We then have $\varepsilon_{\mathscr{C}L}(\drct{J}))\ \sq \clos{J}$ and $\varepsilon_L(\clos{J})=b\in \clos{I}$. Conversely, if such $\clos{J}$ exists, then $\varepsilon_L(\varepsilon_{\mathscr{C}L}(\drct{J}))\kel \varepsilon_L(\clos{J})$ and so $\varepsilon_L(\varepsilon_{\mathscr{C}L}(\drct{J}))\in \clos{I}$.
\end{proof}

\begin{theorem}
\label{thm: C is a comonad}
The triple $(\co,c,\varepsilon)$ forms a comonad on \cat{C}.
\end{theorem}
\begin{proof}
Let us show that the following diagrams commute
$$\xymatrix{
\co\co L \ar[r]^{\mathscr{C} c_L} & \co\co\co L & & \co L & \co\co L\ar[l]_{\varepsilon_{\mathscr{C}L}} \ar[r]^{\mathscr{C}\varepsilon_L} & \co L\\
\co L\ar[u]^{c_L} \ar[r]_{c_L} & \co\co L\ar[u]_{c_{\mathscr{C}L}} & & & \co L \ar[lu]^{1} \ar[u]_{c_L} \ar[ru]_{1} & 
}$$
We have $\varepsilon_{\mathscr{C}L}(c_L(\clos{I}))=\clos{I}$ by Proposition \ref{prop: c is a morphism in C}. For the triangle on the right, we have $\co(\varepsilon_L)c_L=\co(\varepsilon_L)\co(\beta\kappa_L)=\co(\varepsilon_L\beta_L\kappa_L)=\co(1)=1$. This holds since $\rnd(\varepsilon_L)\rnd(\beta\kappa_L)=1$ by Remark \ref{rem: on proxi map to frame map}. For the square on the left
\begin{align*}
c_{\mathscr{C}L}(c_L(\clos{I})) &= c_{\mathscr{C}L}\left(\set{\clos{J}\in\co L\ |\ \varepsilon_L(J)\in \clos{I}}\right)\\
&= \beta_{\mathscr{C}L}\set{\drct{J}\ |\ \varepsilon_L(\varepsilon_{\mathscr{C}L}(\drct{J}))\in \clos{I}}\\
&= \beta_{\mathscr{C}L}\set{\drct{J}\ |\ \varepsilon_{\mathscr{C}L}(\drct{J})\ \sq \clos{J}\text{ for some }\varepsilon_L(\clos{J})\in \clos{I}}
\end{align*}
and on the other hand 
\begin{center}
$\co(c_L)(c_L(\clos{I}))=\beta_{\mathscr{C}L}\set{\drct{J}\ |\ \drct{J}\ \sq c_L(\clos{K})\text{ for some }\varepsilon_L(\clos{K})\in \clos{I}}$.
\end{center}
If $\drct{J}\ \sq c_L(\clos{K})$, then $\varepsilon_{\mathscr{C}L}(\drct{J})\ \sq \varepsilon_{\mathscr{C}L}(c_L(\clos{K}))=\clos{K}$. 
For the converse, assume that $\varepsilon_L(\varepsilon_{\mathscr{C}L}(\drct{J}))\in \clos{I}$. We have $\varepsilon_L(\varepsilon_{\mathscr{C}L}(\drct{J}))\kel b$ for some $b\in \clos{I}$. For all $\clos{K}\in\drct{J}$,\ $\clos{K}\subseteq\varepsilon_{\mathscr{C}L}(\drct{J})$, so that $\varepsilon_L(\clos{K})\leq \varepsilon_L(\varepsilon_{\mathscr{C}L}(\drct{J}))\kel b$ and $\drct{J}\subseteq c_L(\beta_L\kappa_L(b))$. Now if $d\in \clos{I}$ such that $b\kel d$, then $\kappa_L(b)\ll\kappa_L(d)$ and $\beta_L\kappa_L(b)\ \sq\beta_L\kappa_L(d)$. By taking $\clos{K}=\beta_L\kappa_L(d)$, we have $\drct{J}\ \sq c_L(\clos{K})$ and $\varepsilon_L(\clos{K})=b\in \clos{I}$.
\end{proof}
\begin{cor}
\label{cor: submonad}
$(\rnd,r,\varsigma)$ is a submonad\footnote{Horizontal composition of natural transformations is defined in \cite[Section II.3]{Low} and \cite[Chapter II, Section 5]{Mac}.} of $(\co,\varepsilon,c)$. 
\end{cor}
\begin{proof}
We have $c_L\beta_{\mathfrak{R}L}=\beta_{\mathscr{C}L}\rnd(\beta_L\kappa_L)=\beta_{\mathscr{C}L}\rnd(\beta_L)\rnd(\kappa_L)=(\beta_L\circ\beta_L)r_L$ and $\varepsilon_L\beta_L=\varsigma_L$.
\end{proof}

\subsubsection{Coalgebras of the comonad $(\co,\varepsilon,c)$}
Consider two stably compact proximity frames $L$ and $M$. We shall denote by $\alpha_L$ and $\alpha_M$ the maps that are given from Proposition \ref{prop: stably compact proximity frame} without any additional assumption on the nature of the class $\alpha=\set{\alpha_L\ |\ L\in\cat{C}}$ in general.
\begin{prop}
\label{prop: proper map and diagram commutativity}
Consider a frame homomorphism \map{f}{L}{M}.
\begin{enumerate}
\item $\comp{\rnd f}{\alpha_L}=\comp{\alpha_M}{f}$ if and only if $f$ is proper.
\item $\comp{\rnd f}{\alpha_L}=\comp{\alpha_M}{f}$ if and only if $\comp{\co f}{\beta_L\alpha_L}=\comp{\beta_M\alpha_M}{f}$.
\end{enumerate}
\end{prop}
\begin{proof}
\begin{enumerate}
\item Suppose the diagram is commutative and let $a\ll b$ in $L$. Since $\varsigma_L(\alpha_L(a))\in\alpha_L(b)$ we have $\alpha_L(a)\ll \alpha_L(b)$ and $\alpha_M(f(a))\ll \alpha_M(f(b))$. This shows that $f(a)\in\alpha_M(f(b))$. Conversely, suppose that $f$ is proper. By Lemma \ref{lem: rnd join is a nat transformation}, $f\cdot\varsigma_L=\varsigma_M\cdot\rnd f$. By composing with $\alpha_L$ and $\alpha_M$ on the right and on the left respectively, we have $\comp{\rnd f}{\alpha_L}\leq\comp{\alpha_M}{f}$. Now, let $c\in\alpha_M(f(a))$. As $\alpha_M(f(a))$ is round, there is $d\in\alpha_M(f(a))$ with $c\kel d$. Since $f$ preserves joins and $L$ is stably compact, $d\ll\bigvee\set{f(t)\ |\ t\ll a}$ and so there is $t\ll a$ such that $d\leq f(t)\ll f(a)$. Thus $c\in\rnd f(\alpha_L(a))$.
\item This is clear from the fact that $\beta$ is a natural transformation and each component $\beta_M$ is a monomorphism.
\end{enumerate}
\end{proof}
Let us note that the comonad $(\co,\varepsilon,c)$ is a Kock-Z\"oberlein comonad with $\varepsilon_{\mathscr{C}L}\leq\co(\varepsilon_L)$ or equivalently $\co(\varepsilon_L)\dashv\ c_L\dashv\ \varepsilon_{\mathscr{C}L}$\footnote{See \cite[Section II.4.9]{Low} and \cite{Koc}.}. This implies that any coalgebra morphism \map{k}{L}{\co L} is a section that is left adjoint to $\varepsilon_L$.
\begin{prop}
\label{prop: coalgebras of C}
The coalgebras of the comonad $(\co,\varepsilon,c)$ are precisely the proximity frames that are stably compact together with proper frame maps that preserve the proximities.
\end{prop}
\begin{proof}
If $(L,k)$ is a coalgebra, then $L$ is stably compact since $\varepsilon_L\cdot k=1$. Since $k\dashv\ \varepsilon_L$, we have $k=\beta_L\alpha_L$. Proposition \ref{prop: proper map and diagram commutativity} shows that any coalgebra morphism \map{f}{(L,k_L)}{(M,k_M)} is proper. Conversely, if $L$ is stably compact, then \map{\alpha}{L}{\rnd L} exists by Proposition \ref{prop: stably compact proximity frame}, with $\beta_L\alpha\dashv\ \varepsilon_L$. Since $\co(\varepsilon_L)\dashv\ c_L\dashv\ \varepsilon_{\mathscr{C}L}$, $\beta_L\alpha$ is a coalgebra morphism. Finally, if \map{f}{L}{M} is proper in \cat{C}, then \map{f}{(L,\beta_L\alpha_L)}{(M,\beta_M\alpha_M)} is a coalgebra morphism by Proposition \ref{prop: proper map and diagram commutativity}.
\end{proof}

\section{Remark on some familiar categories}
With respect to the category \cat{C} of proximity frames and frame homomorphisms preserving proximity relations, certain categories are of particular interest to us. These are the category of completely regular frames, uniform frames and the category of proximal frames.\\

\paragraph{\textit{Completely regular frames}.} The completely below relation is a proximity which is preserved by frame homomorphisms. This presents the category \cat{CRFrm} formed by such frames as a full subcategory of \cat{C}. Here the restriction of \rnd\ gives the Stone-\v{C}ech compactification. As for the restriction of \co, it is not clear whether \sq\ coincides with the strong inclusion $\ll$ on $\rnd L$. However, the coalgebras of the restriction of \co\ are exactly the compact regular frames.\\

\paragraph{\textit{Proximal frames}} Frames endowed with strong inclusions (\cite{Ban90,Fri}) and frame homomorphisms that preserve them form the category of proximal frames denoted by \cat{ProxFrm}. The embedding $\cat{ProxFrm}\to\cat{C}$ is full and the restriction of \rnd\ coreflects \cat{ProxFrm} onto the category of compact regular frames.\\

\paragraph{\textit{Uniform Frames}} The category \cat{UniFrm} of uniform frames cannot be considered as a subcategory of \cat{C}. There is however a functor \map{G}{\cat{UniFrm}}{\cat{ProxFrm}} (\cite{BanPul,Fri}) such that the composition $\rnd\cdot G$ gives the Samuel compactification (See also \cite[Remark 7.12]{BezHar}).\\

\paragraph{\textit{Frames}} If the order relation $\leq$ on any frame is considered as a proximity, then there is full embedding \map{E}{\cat{Frm}}{\cat{C}}. Here the restriction of \rnd\ coincides with the ideal functor \id. However when \rnd\ coreflects \cat{Frm} onto \cat{StKFrm} inside \cat{C}, it portrays \cat{StKFrm} as a full subcategory. The comonad \co\ with the proximity $\sqsubseteq=\subseteq$, is  on the other hand non-idempotent and the coalgebras  \cat{StKFrm} do not form a full subcategory of \cat{Frm} (See \cite[Theorem 3.3]{BezHar3}). To see this, let $G_{\leq}$ and $G_{\ll}$ be the functors that embed \cat{StKFrm} into \cat{C} as proximity frames with the orders $\leq$ and $\ll$ respectively. We then have the following identities: $\rnd\cdot E=G_{\ll}\cdot\id$ and $\co\cdot E=G_{\leq}\cdot\id$.

That \id\ is comonadic on the category \cat{Frm} was already mentioned in B. Banaschewski and G. C. L. Br\"ummer's paper \cite[Section 3]{BanBru}\footnote{There is a typo in the paragraph of this section. Indeed the co-multiplication is given by $\id(\downarrow)$, not by $\id(\kappa)$ which is an isomorphism.}.

\section*{Acknowledgement}
I am indebted to Themba Dube for having sustained my interest in frames and strong inclusions through various exchanges and conversations over many years. I thank Guram G. Bezhanishvili for having invited me to give an online presentation of this topic at the NMSU Algebra seminar and for alerting me to the paper \cite{Fri} of J. Frith. Finally, I thank the referee for thoughtful suggestions.

\end{document}